\documentclass{amsart}
\usepackage{amssymb}
\usepackage{enumerate}
\usepackage[mathscr]{euscript}
\def\bs{\boldsymbol}
\def\bA{\bs{A}}

\def\bg{\bs{g}}
\def\bk{\bk{g}}
\def\bm{\bs{m}}

\def\bS{\bs{S}}
\def\bbs{\bs{s}}
\def\BD{\bs{D}}
\def\BV{\bs{V}}
\def\bgamma{{\bs \gamma}}
\def\bu{$\bullet$\quad}
\def\Cal{\mathcal}
\def\casesep{, & \text{ if }}
\def\CC{\mathbb C}
\def\CF{\mathcal F}
\def\CK{\mathcal K}
\def\CL{\mathcal L}
\def\CM{\mathcal M}
\def\CO{\mathcal O}
\def\const{\operatorname{const}}
\def\CS{\mathcal S}
\def\DD{\mathbb D}
\def\eps{\varepsilon}
\def\Fh{\mathfrak h}
\newcommand{\myfrac}[2]{\genfrac{}{}{0pt}{0}{#1}{#2}} 
\def\NN{\mathbb N}
\def\ord{\operatorname{ord}}
\def\pder#1#2{\frac{\partial#1}{\partial#2}}
\def\phi{\varphi}
\def\PSH{\mathcal{PSH}}
\def\rho{\varrho}
\def\RR{\mathbb R}

\def\too{\longrightarrow}
\def\tuu{\longmapsto}
\def\TT{\mathbb T}

\def\wdht{\widehat}
\def\ZZ{\mathbb Z}
\makeatletter
\def\th@mytheorem{%
  \let\thm@indent\noindent
  \thm@headfont{\bfseries}
    \itshape
}
\def\th@myremark{%
  \let\thm@indent\noindent
  \thm@headfont{\bfseries}
}
\tagsleft@false
\makeatother
\theoremstyle{mytheorem}
\newtheorem{Theorem}{Theorem}[section]

\theoremstyle{myremark}
\newtheorem{Remark}[Theorem]{Remark}
\newtheorem{Example}[Theorem]{Example}

\numberwithin{equation}{section}

\begin{document}
\title[Remarks on the Sibony functions and pseudometrics]{Remarks on the Sibony functions and pseudometrics} 
\author[M.~Jarnicki]{Marek Jarnicki}
\address{Jagiellonian University, Faculty of Mathematics and Computer Science, Institute of Mathematics,
{\L}ojasiewicza 6, 30-348 Krak\'ow, Poland}
\email{Marek.Jarnicki@im.uj.edu.pl}

\author[P.~Pflug]{Peter Pflug}

\address{Carl von Ossietzky Universit\"at Oldenburg, Institut f\"ur Mathematik,
Postfach 2503, D-26111 Oldenburg, Germany}
\email{Peter.Pflug@uni-oldenburg.de}

\thanks{The research was partially supported by the OPUS grant no. 2015/17/B/ST1/00996 that was financed by the National
Science Centre, Poland}

\begin{abstract}
We discuss some basic properties of the Sibony functions and pseudometrics.
\end{abstract}



\maketitle

\section{Introduction}

Let $G\subset\CC^n$ be a domain. For $a\in G$ let
\begin{align*}
\CM_G(a):&=\{|f|: f\in\CO (G,\DD),\;f(a) =0\},\\ 
\CS^{(p)}_G(a):&=\{\sqrt[p]u:\;\; u:G\too[0,1): \log u\in\PSH(G),\\
&\qquad u\in\Cal C^p(\{a\}),\;\exists_{C>0}: u(z)\leq C\|z-a\|^p,\;z\in G\},\quad p\in\NN,\\
\CK_G(a):&=\{u:\;\;u:G\too [0,1): \log u\in\PSH(G),\\
&\qquad\exists_{C>0}: u(z)\leq C\|z-a\|,\;z\in G\},
\end{align*}
where $\DD\subset\CC$ stands for the unit disc, $\CO(G,\DD)$, resp.~$\PSH(G)$ denote the set of all holomorphic functions on $G$ having values in $\DD$, resp.~the set of all plurisubharmonic functions on $G$, and ``$u\in\Cal C^p(\{a\})$'' means that $u$ is of class $\Cal C^p$ in a neighborhood of $a$  (cf.~\cite{JarPfl2013}, \S\;4.2). Note that $\CS^{(1)}_G(a)$ is different from $\CK_G(a)$ (see Remark \ref{Rem1}\eqref{Rem1bb}). Put 
\begin{multline*}
\CS_G(a):=\CS^{(2)}_G(a)=\{\sqrt{u}:\;\; u:G\too[0,1): \log u\in\PSH(G),\\
u\in\Cal C^2(\{a\}),\;u(0)=0\}.
\end{multline*}
Obviously, $\CM_G(a)\subset\CS_G(a)\subset\CK_G(a)$ and $\CS^{(p)}_G(a)\subset\CK_G(a)$, $p\in\NN$. 
If $\CF\in\{\CM, \CS^{(p)}, \CK\}$, then we define:
\begin{align*}
d_G^\CF(a,z):&=\sup\{v(z):\; v\in\CF_G(a)\},\quad a, z\in G,\\
\delta_G^\CF(a;X):&=\sup\Big\{\limsup_{\lambda\to0}\frac{v(a+\lambda X)}{|\lambda|}: v\in\CF_G(a)\Big\},\quad a\in G,\;X\in\CC^n.
\end{align*}
For $\CF\in\{\CM, \CS, \CK\}$ the families $(d^\CF_G)_G$ and $(\delta^\CF_G)_G$ are \emph{holomorphically contractible}, i.e.

\bu $d^\CF_{\DD}(0,t)=t$, $t\in[0,1)$, \quad $\delta^\CF_{\DD}(0;1)=1$;

\bu for any domains $G\subset\CC^n$, $D\subset\CC^m$ and for any holomorphic mapping $F:G\too D$ we have
\begin{align}
d_D^\CF(F(a),F(z))&\leq d_G^\CF(a,z),\quad a,z\in G, \label{C21}\\
\delta_D^\CF(F(a);F'(a)(X))&\leq\delta_G^\CF(a;X),\quad a\in G,\;X\in\CC^n. \label{C22}
\end{align}

In particular, the families $(d^\CF_G)_G$ and $(\delta^\CF_G)_G$ are invariant under biholomorphic mappings.

If $\CF=\CM$, then we get the \emph{M\"obius pseudodistance} $\bm_G:=d_G^\CM$ and the \emph{Ca\-ra\-th\'e\-odo\-ry--Reiffen pseudometric} 
$\bgamma_G:=\delta_G^\CM$. It is known that
\begin{gather}
\bgamma_G(a,z)=\lim_{\lambda\to0}\frac{\bm_G(a,a+\lambda X)}{|\lambda|}=\max\{|f'(z)(X)|: f\in\CO (G,\DD),\;f(a) =0\}.\label{eq1}
\end{gather}

If $\CF=\CS$, then we get the \emph{Sibony function} $\bbs_G:=d_G^\CS$ and the \emph{Sibony pseudometric} $\bS_G:=\delta_G^\CS$. 
It is known that 
$$
\bS_G(a;X)=\sup\{\sqrt{\CL u(a;X)}: u\in\CS_G(a)\},
$$
where $\CL u(a;X):=\sum_{j,k=1}^n\pder{^2u}{z_j\partial\overline z_k}(a)X_j\overline X_k$ is the \emph{Levi form} (cf.~\cite{JarPfl2013}, Proposition 4.2.16). 
In particular, $\bS_G(a;\cdot)$ is a $\CC$-seminorm.

If $\CF=\CK$, then we get the \emph{pluricomplex Green function} $\bg_G:=d_G^\CK$ and the \emph{Azukawa pseudometric} 
$\bA_G:=\delta_G^\CK$. It is known that $\bg_G(a,\cdot)\in\CK_G(a)$, $\log\bA_G(a;\cdot)\in\PSH(\CC^n)$, and 
\begin{gather}
\bA_G(a;X)=\limsup_{\lambda\to0}\frac{\bg_G(a,a+\lambda X)}{|\lambda|}\quad\text{(cf.~\cite{JarPfl2013}, Lemma 4.2.3).}
\label{eq2}
\end{gather}

If $\CF=\CS^{(p)}$, $p\neq2$, then we get the \emph{higher order Sibony function} $\bbs^{(p)}_G:=d^{\CS^{(p)}}_G$ 
and the \emph{higher order Sibony pseudometric} $\bS^{(p)}_G:=\delta^{\CS^{(p)}}_G$. 

The basic properties of $\bm_G$, $\bgamma_G$, $\bg_G$, 
and $\bA_G$ are well understood. In contrast to that, little is known on properties of $\bS_G$ and almost nothing on $\bbs^{(p)}_G$, $p\in\NN$, and 
$\bS^{(p)}_G$, $p\neq2$.

The main aim of this note is to show that the basic properties of $\bbs^{(p)}_G$ and $\bS^{(p)}_G$ differ essentially from the corresponding properties of $\bm_G$, $\bg_G$, $\bgamma_G$, and $\bA_G$. 

\section{Holomorphic contractibility}

\begin{Remark}\label{Rem1}
\begin{enumerate}[(a)]
\item\label{Rem1a} $\CS^{(p)}_G(a)=\{\sqrt[p]u:\;u:G\too[0,1): \log u\in\PSH(G),\;u\in\Cal C^p(\{a\}),\;\ord_au\geq p\}$, where $\ord_au$ denotes the order of zero of $u$ at $a$.

\item\label{Rem1b} In view of the Taylor formula we have 
$$
\bS^{(p)}_G(a;X)=\sup\big\{(\tfrac1{p!}|u^{(p)}(a)(X)|)^{1/p}: \sqrt[p]{u}\in\CS^{(p)}_G(a)\big\},\quad a\in G,\;X\in\CC^n,
$$
where $u^{(p)}(a):\CC^n\too\RR$ stands for the $p$-th Fr\'echet differential of $u$ at $a$.

\item\label{Rem1bb} In view of \eqref{Rem1b} we get $\bS^{(p)}_G(a;\cdot)\equiv0$ for $p$ odd. In particular, $\bS^{(1)}_\DD(0;1)=0<1=\bA_\DD(0;1)$.

\item\label{Rem1c} $\bbs^{(p)}_G\leq\bg_G$, $\bS^{(p)}_G\leq\bA_G$. In particular, 
$\bbs^{(p)}_{\DD}(0,\lambda)\leq\bg_{\DD}(0,\lambda)=|\lambda|$, $\bS^{(p)}_{\DD}(0;1)\leq\bA_{\DD}(0;1)=1$.

\item\label{Rem1d} If $\bg^{p+\eps}_G(a,\cdot)\in\Cal C^p(\{a\})$ for $0<\eps\ll1$, then $\bg^{1+\eps/p}_G(a,\cdot)\in\CS^{(p)}_G(a)$. Consequently, $\bbs^{(p)}_G(a,\cdot)=\bg_G(a,\cdot)$. 
In particular, $\bbs^{(p)}_{\DD}(0,\lambda)=|\lambda|$, $\lambda\in\DD$. 

\item\label{Rem1e} If $\bg^{2p}_G(a,\cdot)\in\Cal C^{2p}(\{a\})$,  then $\bS^{(2p)}_G(a;\cdot)=\bA_G(a;\cdot)$.
In particular,  $\bS^{(2p)}_\DD(0;1)=1$. 

\item\label{Rem1f} If $F:G\too D$ is holomorphic, then $v\circ F\in\CS^{(p)}_G(a)$ for every $v\in\CS^{(p)}_D(F(a))$. Consequently, 
the family $(\bbs^{(p)}_G)_G$ (resp.~$(\bS^{(p)}_G)_G$) satisfies \eqref{C21} (resp.~\eqref{C22}). 

\item\label{Rem1g} The families $(\bbs^{(p)}_G)_G$ and $(\bS^{(2p)}_G)_G$ are holomorphically contractible. They will be the main objects of our investigation in the sequel.

\item\label{Rem1h} $\bm_G\leq\bbs^{(p)}_G\leq\bg_G$ and $\bgamma_G\leq\bS^{(2p)}_G\leq\bA_G$.
\end{enumerate}
\end{Remark}

\section{Upper semicontinuity}

It is known that for $\CF\in\{\CM, \CK\}$ the functions $G\times G\in(z,w)\tuu d_G^\CF(z,w)$ and $G\times\CC^n\ni(z,X)\tuu\delta_G^\CF(z;X)$
are upper semicontinuous (cf.~\cite{JarPfl2013}, Propositions 2.6.1, 2.7.1(c), 4.2.10(g,k)).
We will prove that in general the functions $\bbs^{(p)}_G(\cdot,z^0)$ and $\bS^{(2p)}_G(\cdot;X^0)$ are not upper semicontinuous
(Examples \ref{Exam1}, \ref{Exam2}). 

Recall that $\bS_G(a;\cdot)$ is a seminorm and therefore it is continuous.
We do not know whether the functions $\bbs_G(a,\cdot)$, $p\in\NN$, and $\bS^{(2p)}_G(a;\cdot)$, $p\geq2$, are upper semicontinuous. 

\begin{Example}[cf.~\cite{JarPfl2013}, Example 4.2.18]\label{Exam1}
Let
$$
G:=\{(z_1,z_2,z_3)\in\CC^3: |z_1|e^{\phi(z_2,z_3)}<1\} 
$$
with
$$
\phi(\xi,\eta):=\sum_{k=1}^\infty\lambda_k\log\Big(\frac{|\xi-a_k|^2+|\eta|}k\Big), \quad(\xi,\eta)\in\CC^2,
$$
where $(a_k)^\infty_{k=1}\subset\DD\setminus\{0\}$ is a dense subset of $\DD$ and $(\lambda_k)_{k=1}^\infty\subset(0,1]$ 
are chosen so that $\phi(0,0)>-\infty$ and $\phi\in\Cal C^\infty(\CC\times\CC_\ast)$, where $\CC_\ast:=\CC\setminus\{0\}$. Note that $G$ is a pseudoconvex Hartogs domain. 

Let $c_t:=(0,0,t)\in G$, $t>0$, $z^0:=(b,0,0)\in G$ with $b\neq0$, and let $X^0:=(1,0,0)$.
We will show that 
\begin{gather*}
\bbs^{(p)}_G(0,z^0)=0<|b|e^{\phi(0,0)}\leq\bbs^{(p)}_G(c_t,z^0),\\
\bS^{(2p)}_G(0;X^0)=0<e^{\phi(0,0)}\leq\bS^{(2p)}_G(c_t;X^0),\quad 0<t\ll1,
\end{gather*}
which shows that \emph{the functions $\bbs^{(p)}_G(\cdot,z^0)$ and $\bS^{(2p)}_G(\cdot;X^0)$ are not upper semicontinuous at $0$.}

Indeed, the function $G\ni(z_1,z_2,z_3)\overset{v}\tuu(|z_1|e^{\phi(z_2,z_3)})^{1+\eps/p}$ belongs to $\CS^{(p)}_G(c_t)$ for all $\eps>0$ and $t>0$. Hence, 
$\bbs^{(p)}_G(c_t,z^0)\geq|b|e^{\phi(0,0)}>0$.
Analogously, the function $G\ni(z_1,z_2,z_3)\overset{v}\tuu|z_1|e^{\phi(z_2,z_3)}$ belongs to $\CS^{(2p)}_G(c_t)$ for all $t>0$.
Hence $\bS^{(2p)}_G(c_t;X^0)\geq\limsup_{\lambda\to0}\frac{v(c_t+\lambda X^0)}{|\lambda|}=e^{\phi(0,t)}\geq e^{\phi(0,0)}>0$.

On the other hand, let $\sqrt[p]{u}\in\CS^{(p)}_G(0)$ (resp.~ $\sqrt[2p]{u}\in\CS^{(2p)}_G(0)$). Since $\CC\times\{a_k\}\times\{0\}\subset G$, we get $u(z_1,a_k,0)=\const(k)$, $z_1\in\CC$, $k\in\NN$. Since $\{0\}\times\CC\times\{0\}\subset G$, 
we get $u(0,z_2,0)=\const=u(0)=0$, $z_2\in\CC$. Thus $u(z_1,a_k,0)=0$, $z_1\in\CC$, $k\in\NN$. 
Since $u\in\Cal C^p(\{0\})$ (resp.~$u\in\Cal C^{2p}(\{0\})$) we conclude that $u=0$ in $U\times\{0\}$, where $U$ is a neighborhood of $(0,0)$. 
Since $\log u\in\PSH(G)$, we get $u(z_1,z_2,0)=0$ for all $(z_1,z_2,0)\in G$. Consequently, $\bbs^{(p)}_G(0,z^0)=0$ (resp.~$\bS^{(2p)}_G(0;X^0)=0$). 
\end{Example}

\begin{Example}
In view of Example \ref{Exam1}, one could expect that perhaps the families $(\bbs^{(p)\ast}_G)_G$ and/or $(\bS^{(2p)\ast}_G)_G$ are holomorphically contractible, where $\bbs^{(p)\ast}_G:=(\bbs^{(p)}_G)^\ast$, $\bS^{(2p)\ast}_G:=(\bS^{(2p)}_G)^\ast$, and ${}^\ast$ denotes the upper semicontinuous regularization. We will prove that unfortunately they are not holomorphically contractible.

Keep the notation from Example \ref{Exam1}. Let 
$$
D:=\{(z_1,z_2)\in\CC^2: (z_1,z_2,0)\in G\},\quad D\ni(z_1,z_2)\overset{F}\tuu(z_1,z_2,0)\in G.
$$ 
Then $\bbs^{(p)\ast}_G(0,z^0)\geq\limsup_{t\to0+}\bbs^{(p)}_G(c_t,z^0)\geq|b|e^{\phi(0,0)}>0$ and 
$\bS^{(2p)\ast}_G(0;X^0)\geq\limsup_{t\to0+}\bS^{(2p)}_G(c_t;X^0)\geq e^{\phi(0,0)}>0$.

On the other hand, let $w^0\in D\cap(\CC\times\DD)$ and let $\sqrt[p]{u}\in\CS^{(p)}_D(\{w^0\})$ (resp.~ $\sqrt[2p]{u}\in\CS^{(2p)}_D(\{w^0\})$).
Since $\CC\times\{a_k\}\subset D$, we get $u(z_1,a_k)=\const(k)$, $z_1\in\CC$, $k\in\NN$. Since $\{0\}\times\CC\subset D$, we get $u(0,z_2)=\const=u(0,0)$, $z_2\in\CC$. Thus $u(z_1,a_k)=\const$, $z_1\in\CC$, $k\in\NN$. 
Since $u\in\Cal C^p(\{w^0\})$ (resp.~$u\in\Cal C^{2p}(\{w^0\})$) we conclude that $u=0$ in $U\times\{0\}$, where $U$ is a neighborhood of $w^0$. Hence, since $\log u\in\PSH(G)$, we get $u(z_1,z_2)=0$ for all $(z_1,z_2)\in D$. Consequently, 
$\bbs^{(p)}_D=0$ on $(D\cap(\CC\times\DD))\times D$ (resp.~$\bS^{(2p)}_D=0$ on $(D\cap(\CC\times\DD))\times\CC^2$). In particular, $\bbs^{(p)\ast}_D(0,(b,0))=0$ (resp.~$\bS^{(2p)\ast}_D(0;(1,0))=0$) and therefore 
\begin{align*}
\bbs^{(p)\ast}_G(F(0,0),F(b,0))&>0=\bbs^{(p)\ast}_D((0,0),(b,0)),\\
\bS^{(2p)\ast}_G(F(0,0);F'(0,0)(1,0))&>0=\bS^{(2p)\ast}_D((0,0);(1,0)).
\end{align*}
\end{Example}

\begin{Example}\label{Exam2}
For $n\geq2$ and $\alpha=(\alpha_1,\dots,\alpha_n)\in\RR^n\setminus\{0\}$ let
$$
\BD_\alpha:=\{z\in\CC^n(\alpha): |z^\alpha|:=|z_1|^{\alpha_1}\cdots|z_n|^{\alpha_n}<1\},
$$
where $\CC^n(\alpha):=\{(z_1,\dots,z_n)\in\CC^n: \forall_{j\in\{1,\dots,n\}}: (\alpha_j<0 \Longrightarrow z_j\neq0)\}$.
Note that $\BD_\alpha$ is a pseudoconvex Reinhardt domain. For $a=(a_1,\dots,a_n)\in\BD_\alpha$ define 
\begin{gather*}
\Xi(a):=\{j\in\{1,\dots,n\}: \alpha_j>0,\;a_j=0\},\\
r(a):=\begin{cases} 1 \casesep \sigma(a)=0 \\ \sum\limits_{j\in\Xi(a)}\alpha_j\casesep \sigma(a)\geq1\end{cases},\quad
\sigma(a):=\#\Xi(a),\quad
\mu(a):=\min\{\alpha_j: j\in\Xi(a)\}.
\end{gather*}
Note that if $\sigma(a)=1$, then $r(a)=\mu(a)$.

The following results are known (cf.~\cite{JarPfl2013}, \S\S\;6.2, 6.3, and \cite{JarPfl2018}, Theorem 1).

\bu If $\alpha_1,\dots,\alpha_n\in\ZZ$ are relatively prime, then
\begin{gather*}
\bm_{\BD_\alpha}(a,z)=\bm_{\DD}(a^\alpha,z^\alpha),\quad
\bg_{\BD_\alpha}(a,z)=\big(\bm_{\DD}(a^\alpha,z^\alpha)\big)^{1/r},\\ 
\bA_{\BD_\alpha}(a;X)=\Big(\bgamma_\DD(a^\alpha;
\tfrac1{r!}\prod_{j\notin\Xi(a)}a_j^{\alpha_j}\cdot\prod_{j\in\Xi(a)}X_j^{\alpha_j})\Big)^{1/r},\quad r=r(a),\\
\bbs_{\BD_\alpha}(a,z)=\begin{cases}\bm_{\DD}(a^\alpha,z^\alpha)\casesep \sigma(a)=0\\
|z^\alpha|^{1/\mu(a)}\casesep \sigma(a)\geq1\end{cases},\\
\bS_{\BD_\alpha}(a;X)=\begin{cases}\bA_{\BD_\alpha}(a;X)\casesep \sigma(a)\leq1\\
0\casesep \sigma(a)\geq2\end{cases},\quad a,z\in\BD_\alpha,\; X\in\CC^n.\\
\end{gather*}

\bu If $\alpha\notin\RR\cdot\ZZ^n$, then
\begin{gather*}
\bm_{\BD_\alpha}\equiv0,\quad 
\bg_{\BD_\alpha}(a,z)=\begin{cases} 0 \casesep \sigma(a)=0\\
|z^\alpha|^{1/r} \casesep \sigma(a)\geq1 \end{cases},\\
\bA_{\BD_\alpha}(a;X)=\begin{cases}0\casesep \sigma(a)=0\\ \big(\prod_{j\notin\Xi(a)}|a_j|^{\alpha_j}\cdot\prod_{j\in\Xi(a)}|X_j|^{\alpha_j}\big)^{1/r}\casesep \sigma(a)\geq1
\end{cases},\quad r=r(a),\\
\bbs_{\BD_\alpha}(a,z)=\begin{cases}0 \casesep \sigma(a)=0\\
|z^\alpha|^{1/\mu(a)}\casesep \sigma(a)\geq1\end{cases},\\
\bS_{\BD_\alpha}(a;X)=\begin{cases}\bA_{\BD_\alpha}(a;X)\casesep \sigma(a)=1\\
0\casesep \sigma(a)\neq1\end{cases},\quad a, z\in\BD_\alpha,\; X\in\CC^n.
\end{gather*}

In particular, if $n=3$ and $\alpha=(1,2,2)$, then
$\bbs_{\BD_\alpha}((0,0,0),z)=|z^\alpha|$ and  $\bbs_{\BD_\alpha}((1/k,0,0),z)=|z^\alpha|^{1/2}$, $k\in\NN$.
Thus, the function $\bbs_{\BD_\alpha}(\cdot, z^0)$ is not upper semicontinuous at $(0,0,0)$ for all $z^0=(z^0_1,z^0_2,z^0_3)\in\BD_\alpha$ 
with $z^0_1z^0_2z^0_3\neq0$.

Notice that using the above effective formulas one may easily construct many other counterexamples.
\end{Example}

\begin{Example}
Keep the notation from Example \ref{Exam2}. Assume that $\alpha_1,\dots,\alpha_n\in\RR_\ast:=\RR\setminus\{0\}$, $a_1\cdots a_s\neq0$, $a_{s+1}=\dots=a_n=0$, 
$s:=n-\sigma(a)$. In particular, $\alpha_{s+1},\dots,\alpha_n>0$.

First observe that if $\sigma(a)\leq1$, then $\bg^{p+\eps}_{\BD_\alpha}(a,\cdot)\in\Cal C^p(\{a\})$ and consequently $\bbs^{(p)}_{\BD_\alpha}(a,\cdot)=\bg_{\BD_\alpha}(a,\cdot)$ 
(Remark \ref{Rem1}\eqref{Rem1d}). Similarly, if $\sigma(a)\leq1$, then $\bg^{2p}_{\BD_\alpha}(a,\cdot)\in\Cal C^\infty(\{a\})$ and consequently 
$\bS^{(2p)}_{\BD_\alpha}(a;\cdot)=\bA_{\BD_\alpha}(a;\cdot)$ (Remark \ref{Rem1}\eqref{Rem1e}). Problems start when $\sigma(a)\geq2$. 
We do not know effective formulas for $\bbs^{(p)}_{\BD_\alpha}(a,\cdot)$, $p\neq2$, and $\bS^{(2p)}_{\BD_\alpha}(a;\cdot)$, $p\geq2$. To illustrate problems we discuss some particular cases. 

\medskip

\bu Assume that $\sigma(a)\geq1$ and $k_j:=\frac{p\alpha_j}{r(a)}\in\NN$, $j=s+1,\dots,n$. Then 
$$
\bg^{2p}_{\BD_\alpha}(a,z)=\prod_{j=1}^s|z_j|^{2p\alpha_j/r(a)}\cdot\prod_{j=s+1}^n|z_j|^{2k_j}
$$ 
and consequently $\bg^{2p}_{\BD_\alpha}(a,\cdot)\in\Cal C^\infty(\{a\})$ which gives
$\bS^{(2p)}_{\BD_\alpha}(a;\cdot)=\bA_{\BD_\alpha}(a;\cdot)$. For example, if $n=2$, $\alpha=(1,1)$, $a=(0,0)$, then 
$\bS^{(4k)}_{\BD_\alpha}(a;X)=|X_1X_2|^{1/2}$, $k\in\NN$.

\medskip

\bu Assume that $\sigma(a)\geq1$ and there exists a $j_0\in\{s+1,\dots,n\}$ such that $\frac{2p\alpha_{j_0}}{r(a)}\notin\NN$.  Then $\bS^{(2p)}_{\BD_\alpha}(a;\cdot)\equiv0$.

Indeed, we may assume that $j_0=n$. Let $r:=r(a)$ and let $k\in\NN_0$ be such that $k<\frac{2p\alpha_n}r<k+1$. In view of Remark \ref{Rem1}\eqref{Rem1b} we have to prove that 
$u^{(2p)}(a)\equiv0$ for all $\sqrt[2p]{u}\in\CS^{(2p)}_{\BD_\alpha}(a)$. Fix such a $u$ and suppose that $u^{(2p)}(a)(X^0)\neq0$ for some $X^0\neq0$. We have
\begin{multline*}
\Big(\frac1{(2p)!}|u^{(2p)}(a)(X)|\Big)^{1/2p}\leq\bS^{(2p)}_{\BD_\alpha}(a;X)\\
\leq\bA_{\BD_\alpha}(a;X)=\Big(\prod_{j=1}^s|a_j|^{\alpha_j}\cdot\prod_{j=s+1}^n|X_j|^{\alpha_j}\Big)^{1/r}.
\end{multline*}
Write $u^{(2p)}(a)(X_1^0,\dots,X_{n-1}^0,tX_n^0)=A_dt^d+\dots+A_0$, $t\in\RR$, with $A_d\neq0$. We have $|A_dt^d+\dots+A_0|\leq\const|t|^{2p\alpha_n/r}$, $t\in\RR$. Taking $t\too\infty$ we get $d\leq k$. On the other hand, taking $t\too0$ we get $A_d=0$; a contradiction.

For example let $n=2$, $\alpha=(q,1)$,  $a=(0,0)$, where
$$
0<q\notin\{\tfrac{2p-k}k: k=1,\dots,2p-1\}\cap\{\tfrac{k}{2p-k}: k=1,\dots,2p-1\}.
$$
Then $\bS^{(2p)}_{\BD_\alpha}(a;\cdot)\equiv0$.

\medskip

\bu As a consequence we conclude that for every $s\in\{0,\dots,n-2\}$ there exists a set $C_s$ dense in $\RR_\ast^s\times\RR_{>0}^{n-s}$ ($\RR_{>0}:=(0,+\infty)$) such that for any $\alpha\in C_s$,
 $a\in\BD_\alpha\cap(\CC_\ast^s\times\{0\}^{n-s})$, and $p\in\NN$ we have $\bS^{(2p)}_{\BD_\alpha}(a;\cdot)\equiv0$.

Indeed, we may put 
$$
C_s:=(\RR_\ast^s\times\RR_{>0}^{n-s})\setminus\bigcup_{\substack{p,k\in\NN:\; k<2p\\ j\in\{s+1,\dots,n\}}}\{\alpha\in\RR^n: 2p\alpha_j=k(\alpha_{s+1}+\dots+\alpha_n)\}.
$$
\end{Example}

Now we turn to discuss a special case where $G\subset\CC^n$ is a complete $n$-circled domain 
(Example \ref{Exam2a}).

\begin{Example}\label{Exam2a} Let $G\subset\CC^n$ be a  complete $n$-circled domain,
i.e.~for any $z=(z_1,\dots,z_n)\in G$ and $\lambda=(\lambda_1,\dots,\lambda_n)\in\overline\DD^n$, the point $\lambda\cdot z:=(\lambda_1z_1,\dots,\lambda_nz_n)$ belongs to
$G$.

(a) Since $\bbs^{(p)}_G(0,\cdot)\leq\bg_G(0,\cdot)$ and the Green function is upper semicontinuous, 
the function $\bbs^{(p)}_G(0,\cdot)$ is continuous at $0$.

\medskip

(b) The function $\bbs^{(p)}_G(0,\cdot)$ is upper semicontinuous in the domain $G\setminus\BV_0$, where 
$\BV_0:=\{(z_1,\dots,z_n)\in\CC^n: z_1\cdots z_n=0\}$.

Indeed, let $M:=\{a\in G: \bbs^{(p)}_G(0,\cdot) \text{ is not upper semicontinuous at }a\}$. 
Since $\bbs^{(p)}_G(0,\cdot)$ is invariant under $n$-rotations (i.e.~under mappings $G\ni z\tuu \lambda\cdot z\in G$,  
$\lambda\in\TT^n$, where $\TT:=\partial\DD$), the set $M$ is also invariant under $n$-rotations.  It is known that $M$ is pluripolar,
i.e.~there exists a $v\in\PSH(\CC^n)$, $v\not\equiv-\infty$, such that $M\subset v^{-1}(-\infty)$ (cf.~\cite{Kli1991}, Theorem 4.7.6).
Suppose that $a=(a_1,\dots,a_n)\in M\setminus\BV_0$. Then $v(\lambda\cdot a)=-\infty$ for all $\lambda\in\TT^n$.
Consequently, by the maximum principle for plurisubharmonic functions,  $v(z_1,\dots,z_n)=-\infty$ for all $|z_j|\leq|a_j|$, $j=1,\dots,n$. Hence $v\equiv-\infty$; a
contradiction.

\medskip

(c) Let $a=(0,\dots,0,a_{s+1},\dots,a_n)=:(0,b)\in G\cap\BV_0$, $1\leq s\leq n-1$, $a_{s+1}\cdots a_n\neq0$.
Define $D:=\{\zeta\in\CC^{n-s}: (\underset{s\times}{\underbrace{0,\dots,0}},\zeta)\in G\}$.
Note that $D$ is a complete $(n-s)$-circled domain 
with $b\in D$. Let $\Fh_D$ denote the Minkowski functional of $D$ ($\Fh_D(\zeta):=\inf\{1/t: t>0,\;t\zeta\in D\}$, $\zeta\in\CC^{n-s}$). Observe that $\Fh_D$ is continuous (because $D$ is $(n-s)$-circled).

\emph{Assume that $\bbs^{(p)}_D(0,b)=\Fh_D(b)$}. Then the function $\bbs^{(p)}_G(0,\cdot)$ is upper semicontinuous at $a$.

Indeed, let $0<R<1$ and $k>0$ be such that $\Fh_D(b)<R$ and $\|b\|<k$. Note that $\{\zeta\in D: \Fh_D(\zeta)<R,\;\|\zeta\|<k\}\subset\subset
D$. Consequently, there exists an $\eps>0$ such that $U:=\{(z',z'')\in\CC^n: \|z'\|<\eps,\;\Fh_D(z'')<R,\;\|z''\|<k\}\subset G$. Then for $z=(z',z'')\in U$ we have
$$
\bbs^{(p)}_G(0,z)\leq\bbs^{(p)}_U(0,z)\leq\bg_U(0,z)\leq\max\Big\{\frac{\|z'\|}\eps, \frac{\Fh_D(z'')}R, \frac{\|z''\|}k\Big\}.
$$
Hence 
$$
\limsup_{z\to a}\bbs^{(p)}_G(0,z)\leq\limsup_{(z',z'')\to(0,b)}
\max\Big\{\frac{\|z'\|}\eps, \frac{\Fh_D(z'')}R, \frac{\|z''\|}k\Big\}=\max\Big\{\frac{\Fh_D(b)}R, \frac{\|b\|}k\Big\}.
$$ 
Letting $R\too1$ and $k\too+\infty$ we get $\limsup_{z\to a}\bbs^{(p)}_G(0,z)\leq\Fh_D(b)$.

On the other side, since the projection $\CC^s\times\CC^{n-s}\ni(z',z'')\tuu z''\in D$ is well-defined, we get
$\Fh_D(b)=\bbs^{(p)}_D(0,b)\leq\bbs^{(p)}_G(0,a)$.

\medskip

(d) Observe that $\bbs^{(p)}_D(0,b)=\Fh_D(b)$ in the case where $D$ is convex. If $s=n-1$, then $D$ is either a disc or the
whole $\CC$. Thus, if $s=n-1$, then  the function $\bbs^{(p)}_G(0,\cdot)$ is upper semicontinuous at each point $a\in\BV_0$ 
of the form $a=(0,\dots,0,a_j,0,\dots,0)\in G$.

\medskip

(e) Consequently, if $n=2$, then the function $\bbs^{(p)}_G(0,\cdot)$ is globally upper semicontinuous.

\medskip

(f) If $G$ is bounded, then the function $\bbs^{(p)}_G(0,\cdot)$ is globally upper semicontinuous.

Indeed, we proceed by induction on $n\geq2$. The case $n=2$ is solved in (e). Suppose the result is true for $n-1\geq2$.  Let $a=(a_1,\dots,a_n)\in G\cap\BV_0$ (see (b)). We may assume that $a_{n-1}\neq0$, $a_n=0$. Define $D:=\{z'\in\CC^{n-1}: (z',0)\in G\}$; $D$ is a bounded complete $(n-1)$-circled domain. Thus, by the inductive assumption, $\bbs^{(p)}_D(0,\cdot)$ is upper semicontinuous. Since $G$ is bounded, for every $0<r<1$ with $a\in rG$ there exists an $\eps>0$ such that $(rD)\times\DD(\eps)\subset\subset G$. Suppose that $G\subset\DD^n(R)$ and let $\eta>0$ be such that 
$rR|\frac{z_n}{z_{n-1}}|<\eps$ for $z\in U:=\{z=(z',z_n)\in a+\DD^n(\eta): z'\in rD\}$. For $z\in U$ consider the holomorphic mapping $F_z:rD\too G$, $F_z(w):=(w,w_{n-1}\frac{z_n}{z_{n-1}})$. We have $\bbs^{(p)}_G(0,F_z(w))\leq\bbs^{(p)}_{rD}(0,w)=\bbs^{(p)}_D(0,w/r)$. In particular,
$\bbs^{(p)}_G(0,z)=\bbs^{(p)}_G(0,F_z(z'))\leq\bbs^{(p)}_D(0,z'/r)$. Thus $\limsup_{z\to a}\bbs^{(p)}_G(0,z)\leq\limsup_{z\to a}\bbs^{(p)}_D(0,z'/r)=\bbs^{(p)}_D(0,a'/r)$. Letting $r\too1-$ (and using once again the upper semicontinuity of $\bbs^{(p)}_D(0,\cdot)$ we get
$\limsup_{z\to a}\bbs^{(p)}_G(0,z)\leq\bbs^{(p)}_D(0,a')\leq\bbs^{(p)}_G(0,a)$ (cf.~(c)).

\medskip

\emph{Note that if $n\geq3$ and $D$ is unbounded, then it is not know whether the function $\bbs^{(p)}_G(0,\cdot)$ is globally upper semicontinuous}.
\end{Example}

\section{Increasing domains property}

Let $(G_k)_{k=1}^\infty$ be sequence of domains in $\CC^n$ such that $G_k\nearrow G$, i.e.~$G_k\subset G_{k+1}$, $k\in\NN$, and 
let $G=\bigcup_{k=1}^\infty G_k$. It is known that if $\CF\in\{\CM, \CK\}$, then $d_{G_k}^\CF\searrow d_G^\CF$ 
and $\delta_{G_k}^\CF\searrow\delta_G^\CF$ (cf.~\cite{JarPfl2013}, Propositions 2.7.1(a), 4.2.10(a)). We will show that 
\emph{this is not true for $\CF=\CS^{(p)}$}.

\begin{Example}\label{Exam3}
Let 
$$
\phi_k(\lambda):=\sum_{s=2}^k\frac1{s^2}\log\Big|\lambda-\frac1s\Big|,\quad k\geq2,\quad 
\phi(\lambda):=\sum_{s=2}^\infty\frac1{s^2}\log\Big|\lambda-\frac1s\Big|,\quad |\lambda|<\frac12.
$$
Observe that $\phi_k\in\PSH$ and $\phi_k\searrow\phi$. Moreover, $\phi_k\in\Cal C^\infty(\tfrac1k\DD)$. Define
\begin{gather*}
G_k:=\{(z_1,z_2)\in\CC^2: |z_1|<1/2,\;|z_2|e^{\phi_k(z_1)}<1\},\\ G:=\{(z_1,z_2)\in\CC^2: |z_1|<1/2,\;|z_2|e^{\phi(z_1)}<1\}.
\end{gather*}
Note that $G_k$ is a Hartogs domain in $\CC^2$, $k\geq2$, and $G_k\nearrow G$. 
For each $k\geq2$ the function 
$G_k\ni(z_1,z_2)\tuu (|z_2|e^{\phi_k(z_1)})^{1+\eps/p}$ belongs to $\CS^{(p)}_{G_k}((0,0))$, $\eps>0$. Hence 
$\bbs^{(p)}_{G_k}((0,0),(0,z_2))\geq|z_2|e^{\phi_k(0)}\geq|z_2|e^{\phi(0)}$ for $|z_2|<e^{-\phi(0)}$.

Analogously, since the function 
$G_k\ni(z_1,z_2)\tuu|z_2|e^{\phi_k(z_1)}$ belongs to $\CS^{(2p)}_{G_k}((0,0))$, we get
$\bS^{(2p)}_{G_k}((0,0);(0,X_2))\geq |X_2|e^{\phi(0)}$ for $X_2\in\CC$ and $k\geq2$. 

Now let $\sqrt[p]{u}\in\CS^{(p)}_G((0,0))$. Since $\{1/s\}\times\CC\subset G$, the Liouville type theorem for subharmonic functions gives
$u(1/s,z_2)=\const(s)=:c_s$, $s\geq2$, $z_2\in\CC$. Since $u(0,0)=0$, we conclude that $c_s\too0$. Since $u$ is continuous near 
$(0,0)$, we get $u(0,z_2)=\lim_{s\to+\infty}u(1/s,z_2)=\lim_{s\to+\infty}c_s=0$, $|z_2|\ll1$. Hence, since $\log u\in\PSH(G)$, we have 
$u(0,z_2)=0$ for all $|z_2|<e^{\phi(0)}$. Consequently, $\bbs^{(p)}_G((0,0),(0,z_2))=0$, $|z_2|<e^{\phi(0)}$, 
and $\bS^{(p)}_G((0,0);(0,X_2))=0$, $X_2\in\CC$.
\end{Example}

\section{Relations between $(\bm_G, \bbs_G, \bg_G)$ and $(\bgamma_G, \bS_G, \bA_G$)}

We will discuss the following two problems. 
\emph{Find a pseudoconvex domain $G\subset\CC^n$, $a\in G$, and $z^0\in G$ (resp.~$X_0\in\CC^n$) such that} 
\begin{gather*}
\bm_G(a,z^0)<\bbs_G(a,z^0)<\bg_G(a, z^0)\\
\text{(resp. } \bgamma_G(a;X_0)<\bS_G(a;X_0)<\bA_G(a;X_0)\text{)}. 
\end{gather*}

\begin{Example}\label{Exam4}
If $\alpha_1,\dots,\alpha_n\in\ZZ$ are relatively prime, $\sigma(a)\geq2$, and $\mu(a)\geq2$, then the domain $G=\BD_\alpha$ 
(cf.~Example \ref{Exam2}) is an example of a pseudoconvex domain (unfortunately, unbounded) 
such that $\bm_G(a,\cdot)<\bbs_G(a,\cdot)<\bg_G(a,\cdot)$ on $\BD_\alpha\setminus\BV_0$.
\emph{It is not know whether there exists a bounded pseudoconvex domain with this property.}
\end{Example}

\begin{Example}\label{Exam5}
Let $G\subset\CC^n$ be a balanced domain (i.e.~$\overline\DD\cdot G=G$) and let $\Fh_G(z)$
be the \emph{Minkowski functional of $G$}. It is known that $G=\{z\in\CC^n: \Fh_G(z)<1\}$. 
Moreover, $\bg_G(0,\cdot)=\Fh_G$ in $G$ $\Longleftrightarrow$ $\bA_G(0;\cdot)\equiv\Fh_G$ $\Longleftrightarrow$
$G$ is pseudoconvex $\Longleftrightarrow$ $\log\Fh_G\in\PSH(\CC^n)$ (cf.~\cite{JarPfl2013}, Proposition 4.2.10(b)).

Let $\wdht G$ be the convex envelope of $G$. It is known that $\wdht G$ is also balanced and
$\Fh_{\wdht G}=\sup\{q: q:\CC^n\too[0,+\infty) \text{ is a $\CC$-seminorm with } q\leq\Fh_G\}$.
Moreover (cf.~\cite{JarPfl2013}, Proposition 2.3.1(d)), $\bgamma_G(0;\cdot)\equiv\Fh_{\wdht G}$.
Thus, if $G$ is pseudoconvex, then $\bgamma_G(0;\cdot)=\Fh_{\wdht G}\geq\bS_G(0;\cdot)$ and hence
$\bgamma_G(0;\cdot)\equiv\bS_G(0;\cdot)\equiv\Fh_{\wdht G}\leq\Fh_G\equiv\bA_G(0;\cdot)$. Consequently, we get the following result.

If $G$ is a balanced pseudoconvex non-convex domain, then
$$
\Fh_{\wdht G}\equiv\bgamma_G(0;\cdot)\equiv\bS_G(0;\cdot)\myfrac{\leq}{\not\equiv}\bA_G(0;\cdot)\equiv\Fh_G.
$$
\end{Example}
In particular, the result solves the problem formulated in Example 4.2.17 from \cite{JarPfl2013}.

\begin{Example}\label{Exam6}
Keep the notation from Example \ref{Exam1}. Then
$$
\bgamma_G(c_t;X_0)<\bS^{(2p)}_G(c_t;X_0)=\bA_G(c_t;X_0)=e^{\phi(0,t)},\quad p\in\NN,\; 0<t\ll1.
$$

Indeed, the function $G\ni(z_1,z_2,z_3)\overset{v}\tuu|z_1|e^{\phi (z_2,z_3)}$ is of the class $\CS^{(2p)}_G(c_t)$, which gives
$$
\bS^{(2p)}_G(c_t;X_0)\geq\limsup_{\lambda\to0}\frac{v(c_t+\lambda X_0)}{|\lambda|}=e^{\phi(0,t)}>0,\quad t>0.
$$
Observe that the mapping $e^{-\phi(0,t)}\DD\ni\lambda\overset{F}\tuu(\lambda,0,t)\in G$ is well-defined. Hence, using the 
holomorphic contractibility, we get
$$
\bA_G(c_t;X_0)=\bA_G(F(0);F'(0)(X_0))\leq\bA_{e^{-\phi(0,t)}\DD}(0;1)=e^{\phi(0,t)}.
$$
Thus $\bS^{(2p)}_G(c_t;X^0)=\bA_G(c_t;X^0)=e^{\phi(0,t)}\geq e^{\phi(0,0)}>0$, $t>0$.

Now, to get the result it suffices to show that $\bgamma_G((0,0,0);X^0)=0$ and then use the continuity of $\bgamma_G(\cdot;X^0)$.
For, let $f\in\CO(G,\DD)$, $f(0,0,0)=0$. Since $\{0\}\times\CC^2\subset G$, the Liouville theorem implies that $f(0,\cdot,\cdot)=\const$.
Since $f(0,0,0)=0$, we get $f(0,\cdot,\cdot)\equiv0$. Since $\CC\times\{a_k\}\times\{0\}\subset G$, we get $f(\cdot,a_k,0)=\const(k)$.
Thus $f(\cdot,a_k,0)\equiv0$. Since the sequence $(a_k)_{k=1}^\infty$ is dense in $\DD$, we conclude that $f=0$ on
$(\CC\times\DD\times\{0\})\cap G$. Thus $f(z_1,0,0)=0$ provided that $|z_1|<e^{-\phi(0,0)}$. Hence $f'(0,0,0)(X^0)=0$ and so 
$\bgamma_G((0,0,0);X^0)=0$.
\end{Example}

\section{Derivative}

Recall that for $\CF\in\{\CM, \CK\}$ we have $\delta^\CF_G(a;X)=\limsup_{\lambda\to0}\frac{d^\CF_G(a,a+\lambda X)}{|\lambda|},\quad a\in G,\;X\in\CC^n$ (cf.~\eqref{eq1} and \eqref{eq2}). It is an open problem whether 
$$
\bS^{(2p)}_G(a;X)=\limsup_{\lambda\to0}\frac{\bbs^{(2p)}_G(a,a+\lambda X)}{|\lambda|},\quad a\in G,\;X\in\CC^n.
$$ 
Observe that 
\begin{multline*}
\bS^{(2p)}_G(a;X)=\sup\Big\{\limsup_{\lambda\to0}\frac{v(a+\lambda X)}{|\lambda|}: v\in\CS^{(2p)}_G(a)\Big\}\\
\leq\limsup_{\lambda\to0}\frac{\bbs^{(2p)}_G(a,a+\lambda X)}{|\lambda|}\leq\limsup_{\lambda\to0}\frac{\bg_G(a,a+\lambda X)}{|\lambda|}=\bA_G(a;X),
\end{multline*}
so the problem is trivial if $\bS^{(2p)}_G(a;X)=\bA_G(a;X)$. 

\bibliographystyle{amsplain}

\end{document}